\documentclass[12pt, leqno]{article}

\def\im{\mathop{\rm Im}}
\def\re{\mathop{\rm Re}}

\newtheorem{theorem}{Theorem}
\newtheorem{lemma}[theorem]{Lemma}
\newtheorem{proposition}[theorem]{Proposition}
\newtheorem{definition}[theorem]{Definition}
\newtheorem{corollary}[theorem]{Corollary}

\newcommand{\begintheorem}{\addtocounter{equation}{1}\begin{theorem}}
\newcommand{\beginlemma}{\addtocounter{equation}{1}\begin{lemma}}
\newcommand{\beginproposition}{\addtocounter{equation}{1}\begin{proposition}}
\newcommand{\begindefinition}{\addtocounter{equation}{1}\begin{definition}}
\newcommand{\begincorollary}{\addtocounter{equation}{1}\begin{corollary}}

\begin{document}

\title{Notes on functions on the unit disk}

\author{Stephen William Semmes \\
	Rice University \\
	Houston, Texas}

\date{}

\maketitle

\tableofcontents

\section{Power series and their derivatives}
\label{section on power series}
\setcounter{equation}{0}

	As usual, ${\bf R}$ denotes the real numbers, ${\bf Z}$
denotes the integers, and ${\bf C}$ denotes the complex numbers.  Thus
each $\alpha \in {\bf C}$ can be expressed as $a + b \, i$, where $a$,
$b$ are real numbers and $i^2 = -1$.  We call $a$, $b$ the real and
imaginary parts of $\alpha$, and denote them $\re \alpha$, $\im
\alpha$, respectively.

	If $z$ is a complex number, with $z = x + y \, i$ where $x$,
$y$ are the real and imaginary parts of $z$, then we write
$\overline{z}$ for the complex conjugate of $z$, which is defined to
be $x - y \, i$.  The complex conjugate of a sum or product of two
complex numbers is equal to the corresponding sum or product of the
complex conjugates of the two complex numbers.  

	The modulus of $z$ is denoted $|z|$ and is the nonnegative
real number defined by $|z|^2 = x^2 + y^2$.  This is equivalent to
saying that $|z|^2 = z \, \overline{z}$.  The triangle inequality for
the modulus states that $|z + w| \le |z| + |w|$ for all complex
numbers $z$, $w$.  We also have that $|z \, w| = |z| \, |w|$ for all
complex numbers $z$, $w$.

	Let us write $\Delta$ for the unit disk in ${\bf C}$, i.e.,
\begin{equation}
	\Delta = \{z \in {\bf C} : |z| < 1 \}.
\end{equation}
The closed unit disk is denoted $\overline{\Delta}$,
\begin{equation}
	\overline{\Delta} = \{z \in {\bf C} : |z| \le 1 \},
\end{equation}
and the unit circle is denoted $\Gamma$,
\begin{equation}
	\Gamma = \{z \in {\bf C} : |z| = 1 \}.
\end{equation}

	We shall be interested in two basic kinds of functions on $\Delta$,
namely, those represented by power series of the form
\begin{equation}
\label{sum_{n = 0}^infty a(n) z^n}
	\sum_{n = 0}^\infty a(n) \, z^n
\end{equation}
for some sequence $\{a(n)\}_{n=0}^\infty$ of complex coefficients,
and those given as
\begin{equation}
\label{sum_{n = 0}^infty a(n) z^n + sum_{n = -infty}^{-1} a(n) overline{z}^n}
	\sum_{n = 0}^\infty a(n) \, z^n 
		+ \sum_{n = -\infty}^{-1} a(n) \overline{z}^{-n}
\end{equation}
for some doubly-infinite sequence $\{a(n)\}_{n=-\infty}^\infty$ of complex
numbers.

	Let us briefly review some of the standard theory of power series
like these.  In order for the series (\ref{sum_{n = 0}^infty a(n) z^n})
to converge for every $z$ in $\Delta$, it is necessary that
\begin{equation}
	|a(n)| \, r^n, \ n \ge 0, \hbox{ be bounded for each } r \in (0,1).
\end{equation}
Similarly, in order for the series in (\ref{sum_{n = 0}^infty a(n) z^n
+ sum_{n = -infty}^{-1} a(n) overline{z}^n}) to converge for all $z$
in $\Delta$, it is necessary that
\begin{equation}
   |a(n)| \, r^n, \ n \in {\bf Z}, \hbox{ be bounded for each } r \in (0,1).
\end{equation}
These conditions are also sufficient for convergence, and indeed they
imply that the corresponding series converge absolutely for each $z$
in $\Delta$, and uniformly on each disk of the form 
\begin{equation}
\label{closed disk { z in {bf C} : |z| le t }}
	\{ z \in {\bf C} : |z| \le t \},
\end{equation}
$0 < t < 1$.  Moreover, the functions on $\Delta$ given by these power
series are continuously differentiable of all orders on $\Delta$.  In
fact, functions of the form (\ref{sum_{n = 0}^infty a(n) z^n}) are
complex analytic on $\Delta$, and functions of the form (\ref{sum_{n =
0}^infty a(n) z^n + sum_{n = -infty}^{-1} a(n) overline{z}^n}) are
harmonic on $\Delta$.

	Let us be more precise.  If we write the complex variable $z$
as $x + i y$, where $x$, $y$ are real, then we can define the
first-order complex derivatives $\partial/\partial z$ and
$\partial/\partial \overline{z}$ in terms of the usual partial
derivatives $\partial/\partial x$ and $\partial/\partial y$ by
\begin{equation}
	\frac{\partial}{\partial z} 
  = \frac{1}{2} \biggl(\frac{\partial}{\partial x} 
			- i \frac{\partial}{\partial y} \biggr), \quad 
	\frac{\partial}{\partial \overline{z}} 
  = \frac{1}{2} \biggl(\frac{\partial}{\partial x} 
			+ i \frac{\partial}{\partial y}\biggr).
\end{equation}
It is easy to see that
\begin{equation}
	\frac{\partial}{\partial z} (z) 
		= \frac{\partial}{\partial \overline{z}} (\overline{z}) = 1,
	\quad \frac{\partial}{\partial z} (\overline{z}) 
		= \frac{\partial}{\partial \overline{z}} (z) = 0.
\end{equation}
The usual Leibnitz rule for first-order derivatives of a product of
functions works for $\partial/\partial z$ and $\partial/\partial
\overline{z}$ just as it does for $\partial/\partial x$ and
$\partial/\partial y$, and thus
\begin{equation}
	\frac{\partial}{\partial z} (z^n) = n \, z^{n-1}, \quad 
		\frac{\partial}{\partial \overline{z}} (\overline{z}^n)
			= n \, \overline{z}^{n-1}, 
	\quad \frac{\partial}{\partial z} (\overline{z}^n) 
		= \frac{\partial}{\partial \overline{z}} (z^n) = 0
\end{equation}
for all positive integers $n$.  Of course $\partial/\partial z$ and
$\partial/\partial \overline{z}$ applied to constant functions is equal
to $0$, which corresponds to the case $n = 0$.

	Notice that
\begin{equation}
	\frac{\partial}{\partial z} \, \frac{\partial}{\partial \overline{z}} 
		= \frac{1}{4} \biggl(\frac{\partial^2}{\partial x^2}
				+ \frac{\partial^2}{\partial y^2} \biggr),
\end{equation}
i.e., $\partial/\partial z \, \partial/\partial \overline{z}$ is equal
to $1/4$ times the usual two-dimensional Laplace operator.  A
continuously differentiable function $f(z)$ defined on some open
subset of ${\bf C}$ is complex-analytic or holomorphic if and only if
$\partial/\partial z f(z) = 0$ on that open set, and similarly a
twice-continuously differentiable function $h(z)$ defined on some open
subset of ${\bf C}$ is harmonic if and only if $(\partial^2 / \partial
x^2 + \partial^2 / \partial y^2) h = 0$ on that open set, which is the
same as saying that $(\partial / \partial z) (\partial / \partial
\overline{z}) h(z) = 0$ on that open set.

	If $f(z)$ is defined on $\Delta$ by a convergent power series
(\ref{sum_{n = 0}^infty a(n) z^n}), then $f(z)$ is holomorphic on $\Delta$
and 
\begin{equation}
	\frac{\partial}{\partial z} f(z) 
		= \sum_{n = 1}^\infty n \, a(n) \, z^{n-1}.
\end{equation}
This follows from standard results about differentiating power series.
Note that the series on the right also converges for all $z$ in $\Delta$.
Similarly, if $h(z)$ is defined on $\Delta$ by convergent series as in
(\ref{sum_{n = 0}^infty a(n) z^n + sum_{n = -infty}^{-1} a(n) overline{z}^n}),
then $h(z)$ is harmonic and
\begin{equation}
	\frac{\partial}{\partial z} h(z) 
		= \sum_{n = 1}^\infty n \, a(n) \, z^{n-1}, \quad
	\frac{\partial}{\partial \overline{z}} h(z)
		= \sum_{n = -\infty}^{-1} (-n) \, a(n) \, \overline{z}^{-n-1}.
\end{equation}
Again these two series converge for all $z$ in $\Delta$.

\section{Continuous functions and the Poisson kernel}
\label{section on continuous functions and the Poisson kernel}
\setcounter{equation}{0}

	Let $h(z)$ be a continuous function on the closed unit disk
$\overline{\Delta}$ which can be represented by a series expansion as
in (\ref{sum_{n = 0}^infty a(n) z^n + sum_{n = -infty}^{-1} a(n)
overline{z}^n}) on the open unit disk $\Delta$, i.e.,
\begin{equation}
\label{series for h(z)}
	h(z) = \sum_{n = 0}^\infty a(n) \, z^n 
		+ \sum_{n = -\infty}^{-1} a(n) \overline{z}^{-n}
			\qquad\hbox{for } z \in \Delta.
\end{equation}
In particular it is assumed that the series in (\ref{series for h(z)})
converges for all $z$ in $\Delta$.

	We would like to determine the coefficients $a(n)$, $n \in
{\bf Z}$, in (\ref{series for h(z)}) from the knowledge of $h$ on the
unit circle $\Gamma$.  Recall that if $\phi(z)$ is a continuous function
on $\Gamma$, then the integral
\begin{equation}
	\int_\Gamma \phi(z) \, |dz|
\end{equation}
is defined, where $|dz|$ is the element of integration by arclength.
We would like to show that
\begin{equation}
\label{formula for a(n), 1}
	a(n) = \frac{1}{2 \pi} \int_\Gamma h(z) \, \overline{z}^n \, |dz|,
\end{equation}
for all integers $n$, which is equivalent to
\begin{equation}
\label{formula for a(n), 2}
	a(n) = \frac{1}{2 \pi} \int_\Gamma h(z) \, z^{-n} \, |dz|,
\end{equation}
since $z^{-1} = \overline{z}$ when $|z| = 1$.  

	Notice first that
\begin{equation}
\label{some integrals on the unit circle}
	\int_\Gamma |dz| = 2 \pi, \ 
	   \int_\Gamma z^m \, |dz| = 0
		\hbox{ when } m \in {\bf Z}, \ m \ne 0.
\end{equation}
If $r \in (0,1)$, then
\begin{equation}
\label{formula for a(n) r^n, 1}
    a(n) \, r^{|n|} 
	= \frac{1}{2 \pi} \int_\Gamma h(r z) \, \overline{z}^n \, |dz|
\end{equation}
for all integers $n$, which is equivalent to
\begin{equation}
\label{formula for a(n) r^{-n}, 2}
	a(n) \, r^{|n|} 
		= \frac{1}{2 \pi} \int_\Gamma h(r z) \, z^{-n} \, |dz|.
\end{equation}
This follows by substituting the series expansion for $h(r z)$ into
the integral and interchanging the order of integration and summation,
which is permissible because of uniform convergence.  One can then
take the limit as $r \to 1$ using the fact that $h(rz)$ tends to
$h(z)$ uniformly for $|z| = 1$ since $h$ is continuous and therefore
uniformly continuous on the closed unit disk, because the latter is
compact.

	If $\zeta$ lies in $\Delta$, then it follows from the series
expansion for $h(\zeta)$ that
\begin{equation}
\label{Poisson integral formula}
	h(\zeta) = \int_\Gamma h(z) \, P(z, \zeta) \, |dz|,
\end{equation}
where
\begin{equation}
\label{Poisson kernel, series expansion}
	P(z, \zeta) 
		= \frac{1}{2 \pi} 
	   \biggl(\sum_{n = 0}^\infty \overline{z}^n \, \zeta^n
	 + \sum_{n = -\infty}^{-1} z^{-n} \, \overline{\zeta}^{-n} \biggr).
\end{equation}
As before, it is permissible to interchange the order of integration
and summation to obtain (\ref{Poisson integral formula}) because of
uniform convergence.  The function $P(z, \zeta)$, $z \in \Gamma$,
$\zeta \in \Delta$, is called the \emph{Poisson kernel}.\index{Poisson
kernel}

	We can rewrite (\ref{Poisson kernel, series expansion}) as
\begin{eqnarray}
\label{Poisson kernel, 2}
	P(z, \zeta) 
	     & = & \frac{1}{2 \pi} 
	  \biggl(-1 + 2 \re \sum_{n=0}^\infty \overline{z}^n \zeta^n \biggr)
								\\
	& = & \frac{1}{2 \pi} 
		\biggl(-1 + 2 \re \frac{1}{1 - \overline{z} \zeta} \biggr).
							\nonumber
\end{eqnarray}
This can be simplified further as
\begin{eqnarray}
\label{Poisson kernel, 3}
	P(z, \zeta) & = & 
  \frac{1}{2 \pi} 
    \biggl(\frac{-|1 - \overline{z} \zeta|^2}{|1 - \overline{z} \zeta|^2}
  + \frac{2 \re (1 - z \overline{\zeta})}{|1 - \overline{z} \zeta|^2} \biggr)
								\\
	& = & \frac{1}{2 \pi} \, 
		\frac{1 - |\zeta|^2}{|z - \zeta|^2}
							\nonumber
\end{eqnarray}
for $z \in \Gamma$, $\zeta \in \Delta$.  

	Notice that
\begin{equation}
	P(z, \zeta) \ge 0
\end{equation}
in particular.  For each $\zeta$ in $\Delta$, we have that
\begin{equation}
\label{integral of P(z, zeta) is equal to 1}
	\int_\Gamma P(z, \zeta) \, |dz| = 1.
\end{equation}
This is an easy consequence of the series expansion (\ref{Poisson
kernel, series expansion}).

	Suppose that we start with an arbitrary real polynomial on
${\bf C}$, which is to say a polynomial in the real variables $x$, $y$
which are the real and imaginary parts of the complex variable $z$.
This is equivalent to saying that we have a general polynomial in $z$
and $\overline{z}$.  The harmonic polynomials are the ones which can
be written as a sum of a polynomial in $z$ alone and a polynomial in
$\overline{z}$ alone.  Harmonic polynomials define functions on the
unit disk of the sort under consideration, with only finitely many
terms in the expansion (\ref{series for h(z)}).

	Given an arbitrary real polynomial on ${\bf C}$, we can
replace it with a harmonic polynomial such that the two polynomials
are equal on the unit circle.  Indeed, to do this one simply replaces
each monomial $z^j \, \overline{z}^k$ with $z^{j-k}$ when $j \ge k$
and with $\overline{z}^{k-j}$ when $j \le k$.  The harmonic polynomial
is uniquely determined by its restriction to the unit circle,
as in the preceding computations, and hence is determined by the
initial real polynomial.
	
	Now suppose that we start with a continuous function $h(z)$ on
$\Gamma$, and we extend $h$ to a function on the closed unit disk
through the formula (\ref{Poisson integral formula}) for $\zeta \in
\Delta$.  By construction, $h$ is represented by the series expansion
(\ref{series for h(z)}) on $\Delta$, with the coefficients $a(n)$ as
in (\ref{formula for a(n), 1}) and (\ref{formula for a(n), 2}).  More
precisely, the $a(n)$'s are bounded in this case, which ensures that
the relevant series converge on $\Delta$.  In fact the extended
function is continuous on all of $\overline{\Delta}$.

	Continuity on $\Delta$ can be viewed as a consequence of the
series expansion, or of continuity properties of $P(z, \zeta)$.  Let
$w$ be an element of $\Gamma$, and let us consider continuity of the
extended function at $w$.  We only need to be concerned about nearby
points $\zeta$ in $\Delta$, because our original function is
continuous on $\Gamma$ by assumption.  In order to deal with these
points $\zeta$, observe that for each $\rho > 0$ we have
\begin{equation}
	\lim_{\zeta \to w \atop \zeta \in \Delta} P(z, \zeta) = 0
\end{equation}
uniformly on $\{\zeta \in \Gamma : |\zeta - w| \ge \rho\}$.  This is
not hard to check, and once one has this, one can also verify the
desired continuity at $w$.

	Of course a continuous function on the unit circle can be
approximated uniformly by a sequence of functions which are
restrictions of real polynomials on the complex plane to the unit
circle.  For each of these approximations we can get a harmonic
extension which is a polynomial, as discussed previously.  The
resulting sequence of harmonic polynomials converges uniformly on the
closed unit disk to a continuous function on the closed unit disk
which is the extension of the initial continuous function on the unit
circle to a continuous function on the closed unit disk which is
harmonic on the open unit disk.

\section{Normal families}
\label{section on normal families}
\index{normal families}
\setcounter{equation}{0}

	Let $\mathcal{H}$ be a collection of functions on $\Delta$
represented by power series as in (\ref{sum_{n = 0}^infty a(n) z^n +
sum_{n = -infty}^{-1} a(n) overline{z}^n}).  We say that $\mathcal{H}$
is a \emph{normal family} if for every $r \in (0,1)$ there is a
real number $M(r)$ such that
\begin{equation}
	|h(z)| \le M(r) \hbox{ when } h \in \mathcal{H} \hbox{ and } |z| \le r.
\end{equation}
Similarly, a collection $\mathcal{C}$ of functions on ${\bf Z}$ is said
to be a normal family if for every $r \in (0,1)$ there is a real
number $C(r)$ such that
\begin{equation}
	|a(n)| \, r^{|n|} \le C(r) \hbox{ when } a \in \mathcal{C}
						\hbox{ and } n \in {\bf Z}.
\end{equation}
	
	If $\mathcal{H}$ is a normal family of functions on $\Delta$
represented by power series as in (\ref{sum_{n = 0}^infty a(n) z^n +
sum_{n = -infty}^{-1} a(n) overline{z}^n}), then from $\mathcal{H}$ we
get a collection of $\mathcal{C}$ of functions on ${\bf Z}$, namely,
the power series coefficients of the functions in $\mathcal{H}$, and
$\mathcal{C}$ is a normal family.  This can be derived from the
formulas (\ref{formula for a(n) r^n, 1}), (\ref{formula for a(n)
r^{-n}, 2}) for the coefficients of a given function.  Conversely, if
$\mathcal{C}$ is a normal family of functions on ${\bf Z}$, then we
get a family $\mathcal{H}$ of functions on $\Delta$ which can be
represented by power series in (\ref{sum_{n = 0}^infty a(n) z^n +
sum_{n = -infty}^{-1} a(n) overline{z}^n}), namely, the series whose
coefficients are in $\mathcal{C}$.  The conditions on the elements of
$\mathcal{C}$ are strong enough to ensure that the corresponding
series converge on all of $\Delta$, and in fact that the family of
functions that results is a normal family.

	Suppose that $\mathcal{H}$ is a normal family of functions on
$\Delta$ represented by power series as in (\ref{sum_{n = 0}^infty
a(n) z^n + sum_{n = -infty}^{-1} a(n) overline{z}^n}), and that
$\mathcal{C}$ is the corresponding normal family of functions on ${\bf
Z}$.  If $\{h_j\}_{j=1}^\infty$ is a sequence of functions in
$\mathcal{H}$ and $h$ is another function in $\mathcal{H}$, then a
natural form of convergence for $\{h_j\}_{j=1}^\infty$ to $h$ is
uniform convergence on every closed disk $\{z \in {\bf C} : |z| \le r
\}$ for $0 < r < 1$.  Similarly, if $\{a_j\}_{j=1}^\infty$ is a
sequence of functions in $\mathcal{C}$, and $a$ is another function in
$\mathcal{C}$, then a natural kind of convergence for
$\{a_j\}_{j=1}^\infty$ to $a$ is convergence at each element of ${\bf
Z}$, i.e., to have $\lim_{j \to \infty} a_j(n) = a(n)$ for all $n$ in
${\bf Z}$.  If $a_j(n)$, $n \in {\bf Z}$, are the power series
coefficients for $h_j$ for each $j$, and if $a(n)$, $n \in {\bf Z}$
are the power series coefficients of $h$, then convergence of
$\{h_j\}_{j=1}^\infty$ to $h$ in the sense just described is
equivalent to convergence of $\{a_j\}_{j=1}^\infty$ to $a$ in the
sense just described for it.  Indeed, to go from convergence of the
$h_j$'s to convergence of the $a_j$'s, one can use the formulae
(\ref{formula for a(n) r^n, 1}), (\ref{formula for a(n) r^{-n}, 2})
for the $a_j(n)$'s and $a(n)$'s in terms of the $h_j$'s and $h$.  For
the converse one can basically just sum the series.  For a fixed $r$,
the normality condition for $\mathcal{C}$ can be used to show that the
total contribution of the coefficients corresponding to large $|n|$ is
small uniformly on the disk $\{z \in {\bf C} : |z| \le r \}$.  Thus
the convergence is largely affected by what happens for finite ranges
of $n$, for which the convergence of the $a_j$'s can be used.

	One also has compactness results for these normal families,
i.e., every sequence in the family has a subsequence which converges
in the sense described in the preceding paragraph.  This is a version
of the well-known Arzela--Ascoli theorem.  Because of the previous
remarks, it is enough to establish this compactness property for
$\mathcal{C}$, where it can be derived using compactness of closed
disks in ${\bf C}$.

\section{Inner products and shift operators}
\label{section on inner products and shift operators}
\setcounter{equation}{0}

	Let us define an inner product on the space of continuous
functions on $\Gamma$ by
\begin{equation}
\label{inner product for functions on unit circle}
	\langle \phi_1, \phi_2 \rangle 
   = \frac{1}{2 \pi} \int_\Gamma \phi_1(z) \, \overline{\phi_2(z)} \, |dz|.
\end{equation}
With respect to this basis, the functions $z^m$, $\overline{z}^n$, and
the constant function $1$ are orthonormal, where $m$ and $n$ run
through the set of positive integers.  These functions are the
eigenfunctions for differentiation on the unit circle.

	As in Section \ref{section on continuous functions and the
Poisson kernel}, for each continuous function $\phi$ on the unit
circle there is a unique continuous function $\Phi$ on the closed
unit disk such that $\Phi$ and $\phi$ are equal on the unit circle,
and $\Phi$ can be represented as 
\begin{equation}
	\Phi(z) = \sum_{n = 0}^\infty a(n, \phi) \, z^n
		   + \sum_{n = -\infty}^{-1} a(n, \phi) \, \overline{z}^{-n}
\end{equation}
for all $z$ in the open unit disk.  For each $r \in (0,1)$, we have that
\begin{equation}
	\frac{1}{2 \pi} \int_\Gamma |\Phi(r z)|^2 \, |dz|
		= \sum_{n = -\infty}^\infty |a(n, \phi)|^2 \, r^{2 |n|}.
\end{equation}
This can be checked by substituting the series expansion for $\Phi$
and computing in a simple way.  Everything converges nicely, because
$r \in (0,1)$.

	We can take the limit as $r \to 1$ to obtain that
\begin{equation}
	\frac{1}{2 \pi} \int_\Gamma |\phi(z)|^2 \, |dz|
		= \sum_{n = -\infty}^\infty |a(n, \phi)|^2.
\end{equation}
In particular, the series on the right side converges, which is to say
that $a(n, \phi)$ lies in $\ell^2({\bf Z})$.  Similarly, for a pair of
continuous functions $\phi_1$, $\phi_2$ on the unit circle, $\Phi_1$,
$\Phi_2$, $a(n, \phi_1)$, and $a(n, \phi_2)$ can be defined just as
for $\phi$, and we have that
\begin{equation}
   \frac{1}{2 \pi} \int_\Gamma \Phi_1(r z) \, \overline{\Phi_2(r z)} \, |dz| 
 = 
 \sum_{n = -\infty}^\infty a(n, \phi_1) \, \overline{a(n, \phi_2)} \, r^{2 |n|}
\end{equation}
for $r \in (0,1)$, and hence
\begin{equation}
   \frac{1}{2 \pi} \int_\Gamma \phi_1(z) \, \overline{\phi_2(z)} \, |dz|
 = \sum_{n = -\infty}^\infty a(n, \phi_1) \, \overline{a(n, \phi_2)}.
\end{equation}

	Consider the linear operator on $\ell^2({\bf Z})$ which sends
a function $a(n)$ in $\ell^2({\bf Z})$ to the function $a(n-1)$.  This
is called the \emph{forward shift operator} on $\ell^2({\bf Z})$.
This operator maps $\ell^2({\bf Z})$ onto itself and preserves the
inner product
\begin{equation}
	\sum_{n = -\infty}^\infty a_1(n) \overline{a_2(n)}
\end{equation}
on $\ell^2({\bf Z})$.  We can use the correspondence between functions
on the unit circle and functions on ${\bf Z}$ to convert this to an
operator acting on functions on the unit circle, and in fact it
corresponds to the operator
\begin{equation}
	\phi(z) \mapsto z \, \phi(z).
\end{equation}
It is easy to check directly that this operator preserves the inner
product of two continuous functions on the unit circle.
	
	Let ${\bf Z}_{+,0}$ denote the set of nonnegative integers,
i.e., ${\bf Z}_{+,0} = {\bf Z}_+ \cup \{0\}$.  Consider the linear
operator on $\ell^2({\bf Z}_{+,0})$ which sends a function $a(n)$ to
the function equal to $a(n-1)$ when $n \ge 1$ and equal to $0$ for $n
= 0$.  Like the previous operator, this one preserves the standard
inner product on $\ell^2({\bf Z}_{+,0})$.  This operator does not map
$\ell^2({\bf Z}_{+,0})$ onto itself, but onto a proper subspace of
codimension $1$.  In terms of the correspondence with functions, one
should now restrict one's attention to functions whose series
expansions on the unit disk are of the form $\sum_{n=0}^\infty a(n) \,
z^n$.  The operator again corresponds to multiplication by $z$ at the
level of functions, and this works on the closed unit disk.

	A bounded linear operator $T$ on a Hilbert space $\mathcal{H}$
is said to be a contraction if its operator norm is less than or equal to
$1$.  If $p(z)$ is a polynomial on the complex plane, so that
\begin{equation}
	p(z) = c_n \, z^n + c_{n-1} \, z^{n-1} + \cdots + c_0,
\end{equation}
where $c_0, \ldots, c_{n-1}, c_n$ are complex numbers, then we can
define a bounded linear operator $p(T)$ on $\mathcal{H}$ by
\begin{equation}
	p(T) = c_n \, T^n + c_{n-1} \, T^{n-1} + \cdots + c_0 \, I,
\end{equation}
where $I$ denotes the identity operator on $\mathcal{H}$.  A
remarkable result of von Neumann and Heinz states that the operator
norm of $p(T)$ is less than or equal to the supremum of $|p(z)|$ over
$z$ in the closed unit disk.  For the shift operators described
before, this can be derived from the representation of the operators
in terms of multiplication by $z$ on functions on the unit circle, in
which case $p$ of the operator corresponds to multiplication by $p(z)$
on the unit circle.

\section{Convex means}
\label{section on convex means}
\setcounter{equation}{0}

	Let $\phi(s)$ be a continuous real-valued function on $[0,\infty)$
which is monotone increasing and convex.  For the purposes of this
section, one may as well assume that $\phi(0) = 0$.

\beginproposition
\label{convex means for functions continuous on overline{Delta}}
Let $h(z)$ be a continuous complex-valued function on the
closed unit disk $\overline{\Delta}$ which can be represented
by a power series 
\begin{equation}
\label{series for h(z) again}
	\sum_{n = 0}^\infty a(n) \, z^n 
		+ \sum_{n = -\infty}^{-1} a(n) \overline{z}^{-n}
\end{equation}
on the open unit disk $\Delta$.  For each real number $r \in (0,1)$
we have that
\begin{equation}
\label{inequality for convex means, 1}
	\frac{1}{2 \pi} \int_\Gamma \phi(|h(r z)|) \, |dz|
		\le \frac{1}{2 \pi} \int_\Gamma \phi(h(z)) \, |dz|.
\end{equation}
Also,
\begin{equation}
\label{inequality for supremum, 1}
	\sup_{z \in \Gamma} |h(r z)| \le \sup_{z \in \Gamma} |h(z)|.
\end{equation}
\end{proposition}

	To prove this, we use the Poisson integral formula
\begin{equation}
\label{Poisson integral formula, repeated}
	h(\zeta) = \int_\Gamma h(z) \, P(z, \zeta) \, |dz|,
\end{equation}
for $\zeta$ in $\Delta$, as in (\ref{Poisson integral formula}).
Because $P(z, \zeta)$ is nonnegative, we have that
\begin{equation}
\label{|h(zeta)| le int_Gamma |h(z)| P(z, zeta) |dz|}
	|h(\zeta)| \le \int_\Gamma |h(z)| \, P(z, \zeta) \, |dz|
\end{equation}
for all $\zeta$ in $\Delta$.  The monotonicity of $\phi(s)$ implies that
\begin{equation}
	\phi(|h(\zeta)|) 
	   \le \phi\biggl(\int_\Gamma |h(z)| \, P(z, \zeta) \, |dz| \biggr).
\end{equation}
Now we apply the convexity of $\phi(s)$ to obtain
\begin{equation}
	\phi(|h(\zeta)|)
		\le \int_\Gamma \phi(|h(z)|) \, P(z, \zeta) \, |dz|.
\end{equation}
This also employs the nonnegativity of $P(z, \zeta)$ and
(\ref{integral of P(z, zeta) is equal to 1}).  In other words, the value
of $\phi$ at an average of numbers is less than or equal to the average
of the values of $\phi$ at those same numbers.

	For $r \in (0,1)$ we obtain that
\begin{equation}
	\frac{1}{2 \pi} \int_\Gamma \phi(|h(r w)|) \, |dw|
  \le \frac{1}{2 \pi} \int_\Gamma \int_\Gamma \phi(|h(z)|) \, P(z, r w) 
							    \, |dz| \, |dw|.
\end{equation}
Using (\ref{integral of P(z, zeta) is equal to 1}) we obtain that
\begin{equation}
	\frac{1}{2 \pi} \int_\Gamma \phi(|h(r w)|) \, |dw|
		\le \frac{1}{2 \pi} \int_\Gamma \phi(|h(z)|) \, |dz|,
\end{equation}
which is exactly what we wanted.

	It is easy to see that
\begin{equation}
	|h(\zeta)| 
	   \le \sup_{a \in \Gamma} |h(a)| \, \int_\Gamma P(z, \zeta) \, |dz|
	   = \sup_{a \in \Gamma} |h(a)|
\end{equation}
for all $\zeta$ in $\Delta$, by (\ref{|h(zeta)| le int_Gamma |h(z)|
P(z, zeta) |dz|}) and (\ref{integral of P(z, zeta) is equal to 1}).
This completes the proof of Proposition \ref{convex means for
functions continuous on overline{Delta}}.

\begincorollary
\label{convex means for functions on Delta}
Let $h(z)$ be a function on the
open unit disk $\Delta$ which can be represented 
by a power series 
\begin{equation}
\label{series for h(z) again, 2}
	\sum_{n = 0}^\infty a(n) \, z^n 
		+ \sum_{n = -\infty}^{-1} a(n) \overline{z}^{-n}
\end{equation}
there.  For each pair of real numbers $r, s \in (0,1)$ with
$r \le s$ we have that
\begin{equation}
\label{inequality for convex means, 2}
	\frac{1}{2 \pi} \int_\Gamma \phi(|h(r z)|) \, |dz|
		\le \frac{1}{2 \pi} \int_\Gamma \phi(h(s z)) \, |dz|.
\end{equation}
Also,
\begin{equation}
\label{inequality for supremum, 2}
	\sup_{z \in \Gamma} |h(r z)| \le \sup_{z \in \Gamma} |h(s z)|.
\end{equation}
\end{corollary}

	This follows from Proposition \ref{convex means for functions
continuous on overline{Delta}} applied to the function $h(s z)$, which
is in fact defined on a neighborhood of the closed unit disk.  If
$h(z)$ is complex analytic function on the open unit disk, which is to
say that $h(z)$ can be represented by a series of the form
\begin{equation}
\label{series in the holomorphic case}
	\sum_{n = 0}^\infty a(n) \, z^n,
\end{equation}
then there are remarkable extensions of these results to larger
classes of monotone functions $\phi(t)$, namely to functions which
can be expressed as a convex function of $\log t$.  In particular,
a basic class of examples of $\phi(t)$ for harmonic functions
are the functions $\phi(t) = t^p$ with $p \ge 1$, while for
holomorphic functions one can allow $p > 0$.


\begin{thebibliography}{ABCde}


\bibitem [Arv] {Arveson} W.~Arveson, {\it A Short Course on Spectral
Theory}, Springer-Verlag, 2002.

\bibitem [AxlBR] {A-B-R} S.~Axler, P.~Bourdon, and W.~Ramey, {\it
Harmonic Function Theory}, second edition, Springer-Verlag, 2001.

\bibitem [Bea1] {Beals1} R.~Beals, {\it Topics in Operator Theory},
University of Chicago Press, 1971.

\bibitem [Bea2] {Beals2} R.~Beals, {\it Advanced Mathematical Analysis:
Periodic Functions and Distributions, Complex Analysis, Laplace
Transform, and Applications}, Springer-Verlag, 1973.

\bibitem [Dou] {Douglas} R.~Douglas, {\it Banach Algebra Techniques in
Operator Theory}, second edition, Springer-Verlag, 1998.

\bibitem [Dur] {Duren} P.~Duren, {\it Theory of $H^p$ Spaces},
Academic Press, 1970.

\bibitem [DurS] {D-S} P.~Duren and A.~Schuster, {\it Bergman Spaces},
American Mathematical Society, 2004.

\bibitem [Gar] {Garnett} J.~Garnett, {\it Bounded Analytic Functions},
Academic Press, 1981.

\bibitem [Gol] {Goldberg} R.~Goldberg, {\it Methods of Real Analysis},
second edition, Wiley, 1976.

\bibitem [GreK] {G-K} R.~Greene and S.~Krantz, {\it Function Theory of
One Complex Variable}, second edition, American Mathematical Society,
2002.

\bibitem [Hal1] {Halmos1} P.~Halmos, {\it A Hilbert Space Problem
Book}, second edition, Springer-Verlag, 1982.

\bibitem [Hal2] {Halmos2} P.~Halmos, {\it Introduction to Hilbert
Space and the Theory of Spectral Multiplicity}, AMS Chelsea
Publishing, 1998.

\bibitem [Hof] {Hoffman} K.~Hoffman, {\it Banach Spaces of Analytic
Functions}, Dover, 1988.

\bibitem [H\"or] {Hormander} L.~H\"ormander, {\it Notions of
Convexity}, Birkh\"auser, 1994.

\bibitem [Kat] {Katznelson} Y.~Katznelson, {\it An Introduction to
Harmonic Analysis}, third edition, Cambridge University Press, 2004.

\bibitem [Koo] {Koosis} P.~Koosis, {\it Introduction to $H_p$ Spaces},
second edition, Cambridge University Press, 1998.

\bibitem [Kra1] {Krantz1} S.~Krantz, {\it Function Theory of Several
Complex Variables}, second edition, AMS Chelsea Publishing, 2002.

\bibitem [Kra2] {Krantz2} S.~Krantz, {\it Complex Analysis: The
Geometric Viewpoint}, second edition, Mathematical Association of
America, 2004.

\bibitem [Nik] {Nikolski} N.~Nikolski, {\it Operators, Functions, and
Systems: An Easy Reading}, volumes 1 and 2, American Mathematical
Society, 2002.

\bibitem [Pel] {Peller} V.~Peller, {\it Hankel Operators and their
Applications}, Springer-Verlag, 2003.

\bibitem [Rud1] {Rudin1} W.~Rudin, {\it Function Theory in 
Polydisks}, Benjamin, 1969.

\bibitem [Rud2] {Rudin2} W.~Rudin, {\it Principles of Mathematical
Analysis}, 3rd edition, McGraw-Hill, 1976.

\bibitem [Rud3] {Rudin3} W.~Rudin, {\it Function Theory in the
Unit Ball of ${\bf C}^n$}, Springer-Verlag, 1980.

\bibitem [Rud4] {Rudin4} W.~Rudin, {\it Real and Complex Analysis},
third edition, McGraw-Hill, 1987.

\bibitem [Rud5] {Rudin5} W.~Rudin, {\it Functional Analysis},
second edition, McGraw-Hill, 1991.

\bibitem [Sar] {Sarason} D.~Sarason, {\it Function Theory on the Unit
Circle}, Virginia Polytecnic Institute and State University, 1978.

\bibitem [SteS] {S-S} E.~Stein and R.~Shakarchi, {\it Complex
Analysis}, Princeton University Press, 2003.

\bibitem [SteW] {S-W} E.~Stein and G.~Weiss, {\it Introduction to
Fourier Analysis on Euclidean Spaces}, Princeton University Press,
1971.

\bibitem [Zyg] {Zygmund} A.~Zygmund, {\it Trigonometric Series},
volumes 1 and 2, third edition, Cambridge University Press, 2002.



\end{thebibliography}
\end{document}